\documentclass[12pt]{amsart}

\usepackage{amsmath,amssymb,bm,amsthm,mathrsfs}
\newtheorem{Def}{'è‹`}[section]

\newtheorem{Them}[Def]{Theorem}
\newtheorem{Lem}[Def]{Lemma}
\newtheorem{Cor}[Def]{Corollary}
\newtheorem{Prop}[Def]{Proposition}
\newtheorem{Examp}[Def]{Example}

\numberwithin{equation}{section}

\newcommand{\Exp}{\operatorname{Exp}}
\newcommand{\rank}{\operatorname{rank}}

\newcommand{\End}{\operatorname{End}}
\newcommand{\ord}{\operatorname{ord}}

\newcommand{\totdeg}{\operatorname{totdeg}}

\newcommand{\Gr}{\operatorname{Gr}}

\newcommand{\NB}{\mathbb{N}}
\newcommand{\AC}{\mathcal{A}}
\newcommand{\FC}{\mathcal{F}}

\newcommand{\DS}{\mathscr{D}}

\title[The Noetherian properties]
{The Noetherian properties of the rings of differential 
operators on central $2$-arrangements}
\author{Norihiro Nakashima}
\date{}
\email{naka\_n@math.sci.hokudai.ac.jp}
\address{Department of Mathematics, Graduate School of Science, 
Hokkaido University, Sapporo, $060$-$0810$, Japan}

\begin{document}

\maketitle

\begin{abstract}
Whereas Holm proved that the ring of differential operators 
on a generic hyperplane arrangement is 
finitely generated as an algebra, 
the problem of its Noetherian properties is still open. 
In this article, after 
proving that the ring of differential operators 
on a central arrangement is right Noetherian 
if and only if it is left Noetherian, 
we prove that the ring of differential operators 
on a central $2$-arrangement is Noetherian. 
In addition, we prove that 
its graded ring associated to the order filtration is not Noetherian 
when the number of the consistuent hyperplanes is greater than $1$.
\vspace{5mm}
\vspace{2mm}
\end{abstract}

\section{Introduction}
Let $K$ be a field of characteristic zero. 
For a commutative $K$-algebra $R$, we inductively define 
$K$-vector spaces of linear differential operators by 
\begin{align*}
\DS^{0}(R)&:=\{\theta\in\End_{K}(R)\mid 
a\in R,\theta a-a\theta=0\},\\
\DS^{m}(R)&:=\left\{ \theta\in\End_{K}(R) \mid 
a\in R,\theta a-a\theta\in\DS^{m-1}(R)\right\}
\quad (m\geq 1).
\end{align*}
We set $\DS(R):=\bigcup_{m\geq 0}\DS^{m}(R)$, and 
we call $\DS(R)$ the ring of differential operators of $R$. 
Let $S:=K[x_{1},\dots,x_{n}]$ denote the polynomial ring. 
It is well known that the ring $\DS(S)$ of 
differential operators of $S$ is the $n$-th Weyl algebra 
$K[x_{1},\dots,x_{n}]\langle 
\partial_{1},\dots,\partial_{n}\rangle$ where 
$\partial_{i}:=\frac{\partial}{\partial x_{i}}$ 
 (see for example \cite{M-R}). 
We use the multi-index notetions, for example, 
$\partial^{\bm{\alpha}}:=\partial_{1}^{\alpha_{1}}
\cdots\partial_{n}^{\alpha_{n}}$ and 
$|\bm{\alpha}|:=\alpha_{1}+\cdots+\alpha_{n}$ 
for $\bm{\alpha}=(\alpha_{1},\dots,\alpha_{n})\in\NB^n$. 
We set $\DS^{(m)}(S):=\bigoplus_{|\bm{\alpha}|=m}
S\partial^{\bm{\alpha}}$ for $m\geq 0$. 
Then the Weyl algebra $\DS(S)$ is decomposed 
into the direct sum of the modules $\DS^{(m)}(S)$ 
of homogeneous differential operators: 
$\DS(S)=\bigoplus_{m\geq 0}\DS^{(m)}(S)$.

There has been a lot of research on finiteness properties of 
the rings of differential operators. 
It is well known that $\DS(R)$ is Noetherian, if 
$R$ is a regular domain (see \cite{M-R}). 
There are some other important classes of algebras 
such that $\DS(R)$ are Noetherian. For example,  
if $R$ is an integral domain of Krull dimension one, 
then $\DS(R)$ is Noetherian (Muhasky \cite{Muh} and 
Smith-Stafford \cite{Sm-St}). 
Saito-Takahashi \cite{Sa-Ta} showed that 
$\DS(R)$ is right Noetherian 
if $R$ is an affine semigroup algebra. 
However, $\DS(R)$ is not Noetherian in general. 
Bernstein-Gel'fand-Gel'fand \cite{BGG} gave an example 
of a ring of differential operators that is neither 
Noetherian nor finitely generated. 

Let $\AC=\left\{ H_{i}\mid i=1,\dots,r \right\}$ 
be a central (hyperplane) arrangement (i.e., every 
hyperplane in $\AC$ contains the origin) in $K^n$. 
Let $I$ be the defining ideal of $\AC$. 
We consider the module $\DS^{(m)}(I)$ of differential operators 
homogeneous of order $m$ that preserve the ideal $I$. 
We call $\DS^{(m)}(I)$ the modules  
of $\AC$-differential operators. 
We find many results about the module $\DS^{(1)}(I)$ 
of $\AC$-derivations in a rich literature 
(see for example \cite{O-T}). 
In contrast, there are only a few literatures about 
the modules of $\AC$-differential operators 
of a higher order. 
Holm \cite{Holm} proved that the 
ring of differential operators of 
the coordinate ring $S/I$ is finitely generated 
when $I$ is the ideal defining 
a generic hyperplane arrangement. 
In this paper, we will prove that $\DS(S/I)$ 
is Noetherian if $n=2$.

In Section \ref{RN-LN}, we prove that 
$\DS(S/I)$ is right Noetherian if and only if 
it is left Noetherian. Thus the Noetherian 
property of $\DS(S/I)$ can be proved by 
the right or left Noetherian property. 

In Section \ref{n=2}, we prove that $\DS(S/I)$ 
is right Noetherian in the case $n=2$. 
Let $R$ be a filtered ring, and $\FC$ the filtration. 
If the graded ring associated to the filtration 
$\FC$ of $R$ is right (left) Noetherian, then 
$R$ is right (left) Noetherian. 
However, the graded ring associated to 
the order filtration of $\DS(S/I)$ is not 
Noetherian if $r\geq 2$ (Example \ref{example}). 
Hence we cannot take this convenient approach to 
prove the Noetherian property of $\DS(S/I)$. 

There is a well-known basis for the module $\DS^{(1)}(I)$ 
of $\AC$-derivation (see for examle \cite{O-T}). 
Holm \cite{Holm-phd} studied the module $\DS^{(m)}(I)$, 
and gave its basis for any order $m$. 
Let $\DS(J)$ denote the subring of $\DS(S)$ consisting of 
the operators preserving an ideal $J$. 
Holm \cite{Holm-phd},\cite{Holm} showed that 
$\DS(I)$ decomposes into the direct sum of $\DS^{(m)}(I)$. 
For an ideal $J$, there is a ring isomorphism:
$$
\DS(S/J)\simeq\DS(J)/J\DS(S)
$$
(see \cite[Theorem 15.5.13]{M-R}). 
Using these facts, we can write any element of $\DS(S/I)$ 
as a linear combination of bases of the modules 
of $\AC$-differential operators. 
This expression is useful to prove Corollary \ref{noeth>r-1}.

We consider a sequence of two-sided ideals 
of $\DS(I)$:
\begin{align*}
I\DS(S)=L_{r}\subseteq L_{r-1} \subseteq \cdots
\subseteq L_{1} \subseteq L_{0}=\DS(I).
\end{align*}
We prove that $\DS(I)/I\DS(S)$ is right Noetherian 
by proving that each $L_{i-1}/L_{i}$ is 
right Noetherian $\DS(I)$-module. 
To show the right Noetherian property 
of $L_{i-1}/L_{i}$, we study a module of lower order operators 
in $L_{i-1}/L_{i}$ and that of higher order operators separately. 

We prove that the right $\DS(I)$-module 
generated by the higher order operators in $L_{i-1}/L_{i}$ 
is Noetherian in Corollary \ref{noeth>r-1}, and that 
the module of lower order operators in $L_{i-1}/L_{i}$ 
is right Noetherian as a right $S$-module 
in Lemma \ref{S-noeth}. In this way, we see that 
$\DS(S/I)$ is Noetherian. 

\section{Differential operators on a central arrangement}
In this section, we fix some notation, 
and we refer to some facts used in Section \ref{n=2}. 
Let $\AC=\left\{ H_{i}\mid i=1,\dots,r \right\}$
 be a central arrangement in $K^n$. 
Fix a polynomial $p_{i}$ defining $H_{i}$, and 
put $Q:=p_{1}\cdots p_{r}$. Thus $Q$ is a product of 
certain homogeneous polynomials of degree $1$. 
We call $Q$ a defining polynomial of $\AC$. 
Let $I$ denote the principal ideal of $S$ generated by $Q$. 

For any ideal $J$ of $S$, we define an $S$-submodule 
$\DS^{(m)}(J)$ of $\DS^{(m)}(S)$ 
and a subring $\DS(J)$ of $\DS(S)$ by 
\begin{align*}
\DS^{(m)}(J)&:=\{\theta\in \DS^{(m)}(S) 
\mid \theta(J) \subseteq J\},\\
\DS(J)&:=\{\theta\in \DS(S)\mid\theta(J)\subseteq J\}. 
\end{align*}
Among others, Holm \cite{Holm} proved 
the following two propositions.
\begin{Prop}[Proposition 4.3 in \cite{Holm}]
\label{p1}
$$
\DS(I)=\bigoplus_{m\geq 0}\DS^{(m)}(I).
$$
\end{Prop}
\begin{Prop}[Proposition 2.4 in \cite{Holm}]
\label{p2}
Suppose that $f_{1},\dots,f_{k}\in S$ are coprime 
to one another. Then
\begin{align*}
\DS(\langle f_{1}\cdots f_{k}\rangle)
=\bigcap_{i=1}^k \DS(\langle f_{i}\rangle).
\end{align*}
\end{Prop}
The following is well known 
(e.g., see \cite[Proposition 2.3]{Holm}).
\begin{Prop}\label{p3}
Let $J$ be the ideal 
of $S$ generated by $f_{1},\dots,f_{k}$,
 and let $\theta\in\DS(S)$ be an operator of 
order $m\geq 1$. Then $\theta\in\DS(J)$ if and only if 
$\theta(x^{\bm{\alpha}}f_{j})\in J$ for 
$|\bm{\alpha}|\leq m-1$ and $j=1,\dots,k$.
\end{Prop}
We use the following lemma in Section \ref{n=2}. 
\begin{Lem}\label{cal}
Let $\delta\in\sum_{i=1}^{n}K\partial_{i}$, and 
let $f_{1},\dots,f_{k}$ be polynomials of degree $1$. 
If $k\leq m$, then 
\begin{align*}
\delta^m f_{1}\dots f_{k}&=\sum^k_{i=0}[m]_{i}
(\sum_{\stackrel{\scriptstyle\Lambda\subseteq \{ 1,\dots ,k\} }
{\sharp \Lambda=i}}\prod_{j\in\Lambda}\delta(f_{j})
\prod_{j\notin \Lambda}f_{j})\delta^{m-i}\\
&=\sum^k_{i=0}[m]_{i}(\frac{1}{i!(k-i)!}
\sum_{\sigma \in S_{k}}\delta(f_{\sigma(1)})\dots 
\delta(f_{\sigma(i)})f_{\sigma(i+1)}\cdots f_{\sigma(k)})
\delta^{m-i},
\end{align*}
where $[m]_{0}:=1$ and 
$[m]_{i}:=m(m-1)\cdots (m-i+1)$ for $i\geq 1$. 
\end{Lem}
\begin{proof}
For any $f\in S$, we see 
$\delta^{\ell} f=f\delta^{\ell}
+{\ell}\delta(f)\delta^{\ell-1}$. 
We can prove the assertion by induction on $k$.
\end{proof}
For a monomial $x^{\bm{\alpha}}\partial^{\bm{\beta}}$ 
in $\DS(S)$, we define its total degree by 
\begin{align}
\totdeg(x^{\bm{\alpha}}\partial^{\bm{\beta}})
=|\bm{\alpha}|-|\bm{\beta}|\label{totdeg}.
\end{align}
For $\theta\in\DS(S)$, we define the total degree 
of $\theta$ as the largest total degree 
of monomials in $\theta$. 
We consider $\DS(S)$ 
a graded ring by the total degree. 

The operator 
$$
\varepsilon_{m}:=\sum_{|\bm{\alpha}|=m}
\frac{m!}{{\bm{\alpha}}!}x^{\bm{\alpha}}\partial^{\bm{\alpha}}
$$ 
is called the Euler operator of order $m$ 
where ${\bm{\alpha}}!=(\alpha_{1}!)\cdots(\alpha_{n}!)$ 
for ${\bm{\alpha}}=(\alpha_{1},\dots,\alpha_{n})$. 
Then $\varepsilon_{1}$ is the Euler derivation, 
and $\varepsilon_{m}=\varepsilon_{1}(\varepsilon_{1}-1)
\cdots(\varepsilon_{1}-m+1)$ 
\cite[Lemma 4.9]{Holm}.
\section{Right Noetherian property 
and left Noetherian property}\label{RN-LN}
Let $Q=p_{1}\cdots p_{r}$ be a 
defining polynomial of a central arrangement 
$\AC$, and let $I=QS$. 
In this section, we will prove that 
the ring $\DS(S/I)$ of differential operators 
is right Noetherian if and only if 
$\DS(S/I)$ is left Noetherian. 
Recall that we have a ring isomorphism
 $\DS(S/I)\simeq \DS(I)/I\DS(S)$ 
(see \cite[Theorem 15.5.13]{M-R}).

Let $0\neq h\in S$, and set $J:=hS$. 
\begin{Lem}\label{DS(J)}
As a ring,
$$
\DS(J)=\DS(S)\cap h\DS(S)h^{-1}.
$$
\end{Lem}
\begin{proof}
Assume that $h\theta h^{-1}\in\DS(S)$ 
with $\theta\in\DS(S)$. For any $f\in S$, 
$$
h\theta h^{-1}(hf)=h\theta (f)\in hS,
$$
which means $h\theta h^{-1}\in\DS(J)$.

Next we will prove the converse inclusion. 
Let $\theta\in\DS(J)$. 
We denote by $K(x_{1},\dots ,x_{n})$ 
the field of fraction of $S$.
Since $h^{-1}\theta h\in K(x_{1},\dots ,x_{n})
\langle\partial_{1},\dots,\partial_{n}\rangle$, 
we can write 
$$
h^{-1}\theta h=\sum_{\bm{\alpha}}
f_{\bm{\alpha}}\partial^{\bm{\alpha}}
$$
with $f_{\bm{\alpha}}
\in K(x_{1},\dots ,x_{n})$. 
We show that $f_{\bm{\alpha}}\in S$ 
for all $\bm{\alpha}$ by induction 
on $|\bm{\alpha}|$.

Since 
$$
f_{0}=h^{-1}\theta h(1)=h^{-1}\theta(h)
\in h^{-1}hS=S,
$$
we have $f_{0}\in S$. 

Assume that $f_{\bm{\alpha}}\in S$ for 
all $\bm{\alpha}$ with $|\bm{\alpha}|<m$. 
For $|\bm{\beta}|=m$, 
$$
h^{-1}\theta h(x^{\bm{\beta}})=
\bm{\beta}!f_{\bm{\beta}}
+\sum_{|\bm{\alpha}|<m}f_{\bm{\alpha}}
\partial^{\bm{\alpha}}(x^{\bm{\beta}}).
$$
Since $\theta\in\DS(J)$, we obtain 
$$
h^{-1}\theta h(x^{\bm{\beta}})=
h^{-1}\theta (hx^{\bm{\beta}})
\in h^{-1}hS =S.
$$
Then $f_{\bm{\beta}}\in S$ by the induction hypothesis. 
Therefore we conclude that 
$h^{-1}\theta h=\sum_{\bm{\alpha}}f_{\bm{\alpha}}
\partial^{\bm{\alpha}}\in\DS(S)$.
\end{proof}
Define an anti-automorphism 
${}^t:\DS(S)\longrightarrow\DS(S)$ 
by ${}^t x_{i}=x_{i},
{}^t\partial_{i}=-\partial_{i}$ 
for $i=1,\dots,n$ 
(we say that ${}^t$ is an anti-automorphism 
if ${}^t$ is an automorphism as a linear map, and 
if ${}^t(\theta \eta) ={}^t\eta{}^t\theta$ 
for any $\theta,\eta\in\DS(S)$). 
It is clear that ${}^t({}^t\theta)=\theta$ 
for any $\theta\in\DS(S)$.

For $\theta\in\DS(J)$, put 
$\theta^{\ast}:=h\,{}^t\theta h^{-1}$. Then 
\begin{align*}
(\DS(J))^{\ast}&=(\DS(S)\cap h\DS(S)
h^{-1})^{\ast}\\
&=h\,{}^t(\DS(S)\cap h\DS(S)h^{-1})h^{-1}\\
&=h\,{}^t\DS(S)h^{-1}\cap{}^t\DS(S)\\
&=h\DS(S)h^{-1}\cap\DS(S)\\
&=\DS(J)
\end{align*}
by Lemma \ref{DS(J)}. Thus
$$
{}^{\ast}:\DS(J)\longrightarrow\DS(J)
$$
is an anti-automorphism. 
If $h\theta\in J\DS(S)$, then
$$
(h\theta)^{\ast}=h\,{}^t (h\theta)h^{-1}=
h\,{}^t\theta hh^{-1}=h\,{}^t\theta\in J\DS(S).
$$
It is clear that $\theta=(\theta^{\ast})^{\ast}$ 
for any $\theta\in\DS(J)$. 
Hence we have $(J\DS(S))^{\ast}=J\DS(S)$. 
Therefore the anti-automorphism ${}^{\ast}$ 
induces an anti-automorphism
$$
{}^{\ast}:\DS(J)/J\DS(S)
\longrightarrow\DS(J)/J\DS(S).
$$
The following is clear from the existence 
of the anti-automorphism ${}^{\ast}$.
\begin{Them}
The ring $\DS(J)/J\DS(S)$ is right Noetherian 
if and only if $\DS(J)/J\DS(S)$ 
is left Noetherian. 
\end{Them}
\begin{Cor}\label{r-l-n}
Let $I$ be the defining ideal of a 
central arrangement. Then 
the ring $\DS(I)/I\DS(S)$ is right Noetherian 
if and only if $\DS(I)/I\DS(S)$ 
is left Noetherian. 
\end{Cor}
%

\section{The case $n=2$}\label{n=2}
In this section, let $n=2$ and $S=K[x,y]$. 
We will prove that the ring $\DS(S/I)\simeq \DS(I)/I\DS(S)$ 
of differential operators is Noetherian. 
We will also prove that, 
in contrast, the graded ring $\Gr\DS(S/I)$ 
associated to the order filtration is not Noetherian 
when $r\geq 2$. 

Put $P_{i}:=\frac{Q}{p_{i}}$ for $i=1,\dots,r$, and define 
\begin{align*}
\delta_{i}:=
\begin{cases}
\partial_{y}&{\rm if}\ p_{i}=ax\quad(a\in K^{\times})\\
\partial_{x}+a_{i}\partial_{y}&{\rm if}\ 
p_{i}=a(y-a_{i}x)\quad(a\in K^{\times}).
\end{cases}
\end{align*}
Then $\delta_{i}(p_{j})=0$ if and only if $i=j$. 
\begin{Prop}[Paper I\hspace{-.1em}I\hspace{-.1em}I,
 Proposition 6.7 in \cite{Holm-phd},
Proposition 4.14 in \cite{Snellman}]
\label{pb-2}
For any $m\geq 1$, $\DS^{(m)}(I)$ is a free 
left $S$-module with basis
\begin{align*}
&\{ \varepsilon_{m}, P_{1}\delta_{1}^m,\dots 
,P_{m}\delta_{m}^m \}\ {\rm if}\ m<r-1,\\
&\{P_{1}\delta_{1}^m,\dots ,P_{r}\delta_{r}^m\}
\ {\rm if}\ m=r-1,\\
&\{P_{1}\delta_{1}^m,\dots ,P_{r}\delta_{r}^m ,
Q\eta_{r+1}^{(m)},\dots ,Q\eta_{m+1}^{(m)} \}\ {\rm if}\ m>r-1,
\end{align*}
where the set $\{\delta_{1}^m,\dots ,\delta_{r}^m ,
\eta_{r+1}^{(m)},\dots ,\eta_{m+1}^{(m)} \}$
forms a $K$-basis for 
$\sum_{|\bm{\alpha}|=m} K\partial^{\bm{\alpha}}$ 
if $m>r-1$.
\end{Prop}
By Proposition \ref{p1}, we have
\begin{align*}
\DS(I)=S&\oplus \bigg{(}\bigoplus_{m=1}^{r-2}
\big{(}S\varepsilon_{m}\oplus SP_{1}\delta_{1}^m
\oplus \cdots \oplus SP_{m}\delta_{m}^m \big{)}\bigg{)}\\
&\oplus \bigg{(} \bigoplus_{m\geq r-1}
\big{(} SP_{1}\delta_{1}^m \oplus \cdots \oplus
SP_{r}\delta_{r}^m \oplus SQ\eta_{r+1}^{(m)}
\oplus \cdots \oplus SQ\eta_{m+1}^{(m)}
\big{)}\bigg{)}.
\end{align*}
For $i=1,\dots,r$, define an additive group
\begin{align*}
L_{i}:=\DS(I)\cap (p_{1}\cdots p_{i})\DS(S).
\end{align*}
\begin{Prop}
For $i=1,\dots,r$, the additive group 
$L_{i}$ is a two-sided ideal of $\DS(I)$.
\end{Prop}
\begin{proof}
It is clear that $L_{i}$ is a right ideal of $\DS(I)$.

To prove that $L_{i}$ is a left ideal of $\DS(I)$, 
by Proposition \ref{p1}, we only need to prove that 
$\DS^{(m)}(I)L_{i}\subseteq L_{i}$ for $m\geq 0$. 
Fix $\theta_{m}\in\DS^{(m)}(I)$. 
For any $j=1,\dots,i$, there exist 
$\eta_{\ell}\in D^{(\ell)}(S)$ such that
\begin{align}
\theta_{m} p_{j}=\eta_{0}+\cdots+\eta_{m}.
\label{th_m-p_j}
\end{align}
We prove that 
$\eta_{\ell}\in p_{j}\bigcap_{i^{\prime}\neq j}
\DS^{(\ell)}(p_{i^{\prime}}S)
\subseteq \DS^{(\ell)}(I)$ for $0\leq \ell\leq m$ 
by induction on $\ell$.

In the case $\ell=0$, let (\ref{th_m-p_j}) act on $1$. Then 
$$
p_{j}S\ni \theta_{m}(p_{j})=\eta_{0}
$$
because $\theta_{m}\in\DS^{(\ell)}(p_{j}S)$ 
by Proposition \ref{p2}. If $\ell>1$, then 
it follows from the induction hypothesis that 
$\eta_{\ell}(x^{\bm{\alpha}})\in p_{j}S$ 
for any $\bm{\alpha}$ since 
$$
p_{j}S\ni \theta_{m}(p_{j}x^{\bm{\alpha}})=
\eta_{0}(x^{\bm{\alpha}})+\cdots+
\eta_{\ell-1}(x^{\bm{\alpha}})+\eta_{\ell}(x^{\bm{\alpha}}).
$$
Therefore 
$\eta_{\ell}\in p_{j}D^{(\ell)}(S)$. Write 
$\eta_{\ell}=p_{j}\eta_{\ell}^{\prime}$. For any $i^{\prime}\neq j$ 
and $|\bm{\alpha}|=\ell-1$, 
it also follows from the induction hypothesis that 
$p_{j}\eta_{\ell}^{\prime}(p_{i^{\prime}}x^{\bm{\alpha}})=
\eta_{\ell}(p_{i^{\prime}}x^{\bm{\alpha}})\in p_{i^{\prime}}S$ 
since 
$$
p_{i^{\prime}}S\ni 
\theta_{m}(p_{j}p_{i^{\prime}}x^{\bm{\alpha}})=
\eta_{0}(p_{i^{\prime}}x^{\bm{\alpha}})
+\cdots+\eta_{\ell-1}(p_{i^{\prime}}x^{\bm{\alpha}})
+\eta_{\ell}(p_{i^{\prime}}x^{\bm{\alpha}}).
$$
Since $p_{j}$ and $p_{i^{\prime}}$ are coprime, we see that 
$\eta_{\ell}^{\prime}(p_{i^{\prime}}x^{\bm{\alpha}})\in p_{i^{\prime}}S$. 
So $\eta_{\ell}^{\prime}\in \DS(p_{i^{\prime}}S)$ 
by Proposition \ref{p3}, and 
$\eta_{\ell}\in p_{j}\bigcap_{i^{\prime}\neq j}
\DS^{(\ell)}(p_{i^{\prime}}S)$. Thus 
$\theta_{m}p_{j}\in p_{j}\bigcap_{i^{\prime}\neq j}\DS(p_{i^{\prime}}S)$. 
Then we conclude that 
$$
\DS(I)p_{1}\cdots p_{i}\DS(S)\subseteq p_{1}\cdots p_{i}\DS(S).
$$
\end{proof}
By Proposition \ref{p1}, $L_{i}$ is decomposed as follows:
\begin{align*}
L_{i}=\bigoplus_{m\geq 0}L^{(m)}_{i},
\end{align*}
where $L^{(m)}_{i}:=\DS^{(m)}(I)\cap(p_{1}\cdots p_{i})\DS^{(m)}(S)$. 
We consider a sequence 
\begin{align}
I\DS(S)=L_{r}\subseteq L_{r-1} \subseteq \cdots
\subseteq L_{1} \subseteq L_{0}=\DS(I) \label{seqL}
\end{align}
of two-sided ideals of $\DS(I)$. 
If a right $\DS(I)$-module $L_{i-1}/L_{i}$ 
is Noetherian for any $i$, 
then $\DS(I)/I\DS(S)$ is a right Noetherian ring. 
Now we fix $i$, and we will prove that 
$L_{i-1}/L_{i}$ is right Noetherian.

As an $S$-module, 
\begin{align*}
L_{i-1}/L_{i}=\bigoplus_{m\geq 0}
\big{(}L_{i-1}^{(m)}+L_{i}/L_{i}\big{)}
\simeq\bigoplus_{m\geq 0}
\big{(}L_{i-1}^{(m)}/L_{i}^{(m)}\big{)}.
\end{align*}
Put
\begin{align*}
\big{(}L_{i-1}/L_{i}\big{)}^{< r-1}&:=
\bigoplus_{m< r-1}
\big{(}L_{i-1}^{(m)}/L_{i}^{(m)}\big{)},\\
\big{(}L_{i-1}/L_{i}\big{)}^{\geq r-1}&:=
\bigoplus_{m\geq r-1}
\big{(}L_{i-1}^{(m)}/L_{i}^{(m)}\big{)}.
\end{align*}
Then $L_{i-1}/L_{i}$ is decomposed as an $S$-module:
\begin{align}
L_{i-1}/L_{i}=\big{(}L_{i-1}/L_{i}\big{)}^{< r-1}
\oplus \big{(}L_{i-1}/L_{i}\big{)}^{\geq r-1}.
\end{align}
We will study $\big{(}L_{i-1}/L_{i}\big{)}^{< r-1}$ 
and $\big{(}L_{i-1}/L_{i}\big{)}^{\geq r-1}$ separately.

First we argue the part of order $\geq r-1$.
\begin{Lem}\label{lemLi}
Assume that $m\geq r-1$. As a left $S$-module, 
\begin{align*}
L_{i}^{(m)}=SQ\delta_{1}^m \oplus \cdots \oplus
SQ\delta_{i}^m \oplus 
&SP_{i+1}\delta_{i+1}^m \oplus \cdots \oplus
SP_{r}\delta_{r}^m \\
&\oplus SQ\eta_{r+1}^{(m)}
\oplus \cdots \oplus SQ\eta_{m+1}^{(m)}.
\end{align*}
\end{Lem}
\begin{proof}
Recall that $P_{i}=\frac{Q}{p_{i}}$. 
We see the assertion by Proposition \ref{pb-2} 
and the definition of $L_{i}$.
\end{proof}
\begin{Prop}\label{lemi-1/i}
As a left $S$-module,
\begin{align*}
\big{(}L_{i-1}/L_{i}\big{)}^{\geq r-1}
=\bigoplus_{m\geq r-1}
\big{(}SP_{i}\delta_{i}^m+L_{i}^{(m)}/
L_{i}^{(m)}\big{)}
\simeq\bigoplus_{m\geq r-1}
\big{(}SP_{i}\delta_{i}^m / SQ\delta_{i}^m\big{)}.
\end{align*}
\end{Prop}
\begin{proof}
By Lemma \ref{lemLi}, 
$L_{i-1}^{(m)}=SP_{i}\delta_{i}^m +L_{i}^{(m)}$ 
for $m\geq r-1$. Then as a left $S$-module
\begin{align*}
\big{(}L_{i-1}/L_{i}\big{)}^{\geq r-1}
=\bigoplus_{m\geq r-1}
\big{(}SP_{i}\delta_{i}^m+L_{i}^{(m)}/
L_{i}^{(m)}\big{)}
\simeq\bigoplus_{m\geq r-1}
\big{(}SP_{i}\delta_{i}^m/
L_{i}^{(m)}\cap SP_{i}\delta_{i}^m\big{)}.
\end{align*}
It remains to prove that 
\begin{align*}
L_{i}^{(m)}\cap SP_{i}\delta_{i}^m
=SQ\delta_{i}^m \subseteq L_{i-1}^{(m)}
\end{align*}
for $m\geq 1$. It is clear that
$SQ\delta_{i}^m\subseteq
L_{i}^{(m)}\cap SP_{i}\delta_{i}^m$. 
Conversely, suppose that 
$fP_{i}\delta_{i}^m\in L_{i}^{(m)}$ with $f\in S$. 
Then $fP_{i}\delta_{i}^m\in p_{1}\cdots p_{i}D^{(m)}(S)$. 
Since the polynomials $p_{i},\dots ,p_{r}$ are 
coprime to one another, we have $f\in p_{i}S$. Thus 
$L_{i}^{(m)}\cap SP_{i}\delta_{i}^m
\subseteq SQ\delta_{i}^m$.
\end{proof}
We define a left $S$-module
\begin{align}
E_{i}:=\bigoplus_{m\geq 0}\big{(}
SP_{i}\delta_{i}^m+L_{i}^{(m)}/L_{i}^{(m)}\big{)}\simeq
\bigoplus_{m\geq 0}\big{(}SP_{i}\delta_{i}^m
/SQ\delta_{i}^m\big{)}.\label{E_i}
\end{align}
Note that $(SP_{i}+L_{i}^{(0)})/L_{i}^{(0)}
\simeq SP_{i}/SP_{i}\cap L_{i}^{(0)}=SP_{i}/SQ$. 
By Proposition \ref{lemi-1/i}, we may identify 
$\big{(}L_{i-1}/L_{i}\big{)}^{\geq r-1}$ with 
the $S$-submodule of $E_{i}$ of order $m\geq r-1$. 
For $g\in S$, we have
\begin{align*}
P_{i}\delta_{i}^m(Qg)
=Q\delta_{i}^m(\frac{Q}{p_{i}}g)\in QS,
\end{align*}
and hence $P_{i}\delta_{i}^m \in \DS^{(m)}(I)$. 
Hence $E_{i}$ is a left $S$-submodule of $L_{i-1}/L_{i}$. 
Moreover, the following proposition is true:
\begin{Prop}\label{P.Ei.r}
The module $E_{i}$ is a right 
$\DS(I)$-submodule of $L_{i-1}/L_{i}$.
\end{Prop}
\begin{proof}
We only need to check the right multiplication by 
the elements of $S$ and the 
bases for $\DS(I)$ in Proposition \ref{pb-2}.

Let $m\geq 1$. For $g\in S$, we have 
$$
\delta_{i}^m \cdot g\in S+\sum_{\ell=1}^m S\delta_{i}^{\ell},
$$
and hence $P_{i}\delta_{i}^m \cdot S \subseteq E_{i}$. 

We show that $E_{i}$ is closed under the right action 
of the elements of bases for $\DS(I)$. We only need to check 
the right multiplication by 
the elements $P_{i}\delta_{i}^{\ell},P_{j}\delta_{j}^m(j\neq i),
\varepsilon_{\ell},Q\eta_{j}^{(\ell)}$. 
For $m\geq 1$, we have 
\begin{align*}
P_{i}\delta_{i}^m \cdot P_{i}\delta_{i}^{\ell}
&=P_{i}(\delta_{i}^m \cdot P_{i})\delta_{i}^{\ell}
\in\bigoplus_{m\geq 0}\big{(}
SP_{i}\delta_{i}^m+L_{i}^{(m)}\big{)},\\
P_{i}\delta_{i}^m \cdot P_{j}\delta_{j}^{\ell}
&=Q\delta_{i}^m \cdot \frac{P_{j}}{p_{i}}\delta_{j}^l
\in\DS(I)\cap (p_{1}\cdots p_{i})\DS(S)=L_{i},\\
P_{i}\delta_{i}^m \cdot Q\eta_{j}^{(\ell)}
&=Q\delta_{i}^m \cdot \frac{Q}{p_{i}}\eta_{j}^{(\ell)}
\in\DS(I)\cap (p_{1}\cdots p_{i})\DS(S)=L_{i}.
\end{align*}
It remains to show that $E_{i}$ is closed under the 
right multiplication by $\varepsilon_{\ell}=\varepsilon_{1}
(\varepsilon_{1}-1)\cdots (\varepsilon_{1}-\ell+1)$. 
We consider the Euler derivation $\varepsilon_{1}$. 
We may assume $p_{i}=y-ax\ (a\in K^{\times})$. 
Recall $\delta_{i}=a^{-1}\partial_{x} +\partial_{y}$. Since
\begin{align*}
\varepsilon_{1}&=x\partial_{x}+y\partial_{y}\\
&=a^{-1}ax\partial_{x}+y\partial_{y}
+a^{-1}y\partial_{x}-a^{-1}y\partial_{x}\\
&=a^{-1}(ax-y)\partial_{x}+y(\partial_{y}+a^{-1}\partial_{x})\\
&=-a^{-1}p_{i}\partial_{x}+y\delta_{i},
\end{align*}
we have, for any $m\geq 0$,
\begin{align*}
P_{i}\delta_{i}^m \cdot \varepsilon_{1}
=P_{i}\delta_{i}^m \cdot
(-a^{-1}p_{i}\partial_{x}+y\delta_{i})
=-a^{-1}Q\delta_{i}^m \partial_{x}
+yP_{i}\delta_{i}^{m+1}
+mP_{i}\delta_{i}^m.
\end{align*}
We see that $-a^{-1}Q\delta_{i}^m \partial_{x} \in L_{i}$, 
and that the remaining terms belong to $SP_{i}\delta_{i}^{m+1}$ 
and $SP_{i}\delta_{i}^m$, respectively. It follows that 
$$
P_{i}\delta_{i}^m \cdot \varepsilon_{\ell}
\in\bigoplus_{m\geq 0}\big{(}
SP_{i}\delta_{i}^m+L_{i}^{(m)}\big{)}.
$$
Hence $E_{i}\cdot\varepsilon_{\ell}\subseteq E_{i}$. 
This completes the assertion.
\end{proof}
As a left $S$-module,
$$
\big{(}L_{i-1}/L_{i}\big{)}^{\geq r-1}
\subseteq E_{i}.
$$
The right $\DS(I)$-module generated by 
$\big{(}L_{i-1}/L_{i}\big{)}^{\geq r-1}$ 
is a $\DS(I)$-submodule of $E_{i}$ 
by Proposition \ref{P.Ei.r}:
\begin{align*}
\big{(}L_{i-1}/L_{i}\big{)}^{\geq r-1}
\cdot \DS(I) \subseteq E_{i}.
\end{align*}
If we prove that $E_{i}$ is a right Noetherian 
$\DS(I)$-module, then 
$\big{(}L_{i-1}/L_{i}\big{)}^{\geq r-1}\cdot \DS(I)$ 
is Noetherian as a $\DS(I)$-module. 
We will prove that $E_{i}$ is a right Noetherian.

We define a left action of $S/p_{i}S$ on $E_{i}$ by 
$$
\overline{f}\cdot \overline{\theta}
=\overline{f\theta}
$$
for $\overline{f}\in S/p_{i}S$ and $\overline{\theta}\in E_{i}$.
This is well-defined, since
\begin{align*}
f\theta-g\theta^{\prime}
=\frac{(f-g)(\theta+\theta^{\prime})}{2}
+\frac{(f+g)(\theta-\theta^{\prime})}{2}
\in L_{i}
\end{align*}
for $f,g\in S$ and $\theta,\theta^{\prime}\in
\bigoplus_{m\geq 0}\big{(}SP_{i}\delta_{i}^m+L_{i}^{(m)}\big{)}$ 
with $f-g\in p_{i}S$ and 
$\theta-\theta^{\prime}\in L_{i}$. 
Thus $E_{i}$ is a left $S/p_{i}S$-module. 
We may assume that $p_{i}=y-ax$ with $a\neq 0$. 
Then $E_{i}$ is a $K$-vector space with a basis 
$\left\{\overline{y}^{\alpha}\cdot
\overline{P_{i}\delta_{i}^m}
\mid\alpha\in\NB, m\geq 0\right\}$. 

Define an exponent by
\begin{align*}
\exp(\overline{y}^{\alpha}\cdot
\overline{P_{i}\delta_{i}^m}):=(\alpha+r-1,m)
\end{align*}
for an element of the basis above. 
We call $\overline{y}^{\alpha}\cdot
\overline{P_{i}\delta_{i}^m}$ a monomial of $E_{i}$. 
Let $\overline{\theta_{1}}$ and $\overline{\theta_{2}}$ be 
two monomials of $E_{i}$ with 
$\exp(\overline{\theta_{1}})=(\alpha_{1},m_{1})$ and 
$\exp(\overline{\theta_{2}})=(\alpha_{2},m_{2})$. 
We define a total order in the set of exponents of monomials by 
$$
\exp(\overline{\theta_{1}})<\exp(\overline{\theta_{2}}),
$$
if $m_{1}<m_{2}$, or if $m_{1}=m_{2}$ and $\alpha_{1}<\alpha_{2}$. 
For $\overline{\theta} \in E_{i}$, write $\overline{\theta}$ 
as a linear combination of monomials. 
Then we define an exponent of $\overline{\theta}$ as 
the largest exponent of a monomial in 
$\overline{\theta}$ with a nonzero coefficient, 
and we denote it by $\exp(\overline{\theta})$. 
For a subset $X$ of $E_{i}$, set
\begin{align*}
\Exp(X):=\left\{ \exp(\overline{\theta})
\mid\overline{\theta}\in X\right\}.
\end{align*}

Throughout the remaining of this section, we 
write $\theta\in E_{i}$ instead of $\overline{\theta}$ 
for simplicity. 

\begin{Lem}\label{lemexp}
Let $M_{1},M_{2}$ be right $\DS(I)$-submodules 
of $E_{i}$. If $M_{1}\subseteq M_{2}$ and 
$\Exp(M_{1})=\Exp(M_{2})$, then
$$
M_{1}=M_{2}.
$$
\end{Lem}
\begin{proof}
Suppose that $M_{1}\subsetneq M_{2}$. 
We can take an element $\theta \in M_{2}\setminus M_{1}$ 
such that $\exp(\theta)$ is the smallest exponent
 in $M_{2}\setminus M_{1}$. 

Since $\exp(\theta)\in\Exp(M_{2})=\Exp(M_{1})$, 
there exists $\eta \in M_{1}$ such that 
$\exp(\eta)=\exp(\theta)$. Then 
$$
\exp(\theta-c\eta)<\exp(\theta)
$$
for some $c\in K^{\times}$. 
We have $\theta-c\eta \in M_{2}\setminus M_{1}$ 
since $\theta\not\in M_{1}$. 
This is a contradiction to the minimality. 
\end{proof}
\begin{Lem}\label{lem<exp}
Let $M\neq 0$ be a right $\DS(I)$-submodule of $E_{i}$. 
If $(k,m)\in \Exp(M)$, then 
\begin{align*}
\left\{(k+a,m),(k+b,m+m^{\prime})\mid
a\geq 0,b\geq r-1,m^{\prime}\geq 1\right\}
\subset\Exp(M).
\end{align*}
\end{Lem}
\begin{proof}
By the assumption, there exists $\theta\in M$ 
such that $\exp(\theta)=(k,m)$. Write 
$\theta=fP_{i}\delta_{i}^m+\theta^{\prime}$ 
with $\exp(\theta^{\prime})<\exp
\left(fP_{i}\delta_{i}^m\right)$. 
The multiplication $\theta\cdot y^a$ 
belongs to $M$, since $S\subseteq\DS(I)$. 
Thus we see that $(k+a,m)\in\Exp(M)$ for all $a\geq 0$.

Fix $1\leq j\neq i\leq r$, $b\geq r+1$, 
and $m^{\prime}\geq 1$. We can write
$$
fP_{i}\delta_{i}^m\cdot 
p_{j}^{b-r+1} P_{i}\delta_{i}^{m^{\prime}}
=f(p_{j}^{b-r+1} P_{i}) P_{i}\delta_{i}^{m+m^{\prime}}
+\eta.
$$
for some $\eta\in E_{i}$ with 
$\exp(\eta)<\exp\left(
f(p_{j}^{b-r+1} P_{i}) P_{i}\delta_{i}^{m+m^{\prime}}\right)$. 
Since $p_{j}^{b-r+1} P_{i}\not\in p_{i}S$, 
we see that $\exp(\theta\cdot p_{j}^{b-r+1}
 P_{i}\delta_{i}^{m^{\prime}})=(k+b,m+m^{\prime})$.
Therefore $(k+b,m+m^{\prime})\in\Exp(M)$ since 
$\theta\cdot p_{j}^{b-r+1}P_{i}\delta_{i}^{m^{\prime}}\in M$.
\end{proof}
Now we induce the total degree (\ref{totdeg}) of $\DS(S)$ 
to those of $\DS(I)$ and $E_{i}$. 
Then $E_{i}$ becomes a graded $\DS(I)$-module 
by the total degree. For monomials of $E_{i}$, 
we denote the total degree by 
\begin{align*}
\totdeg(y^{\alpha})=\alpha,\,
\totdeg(y^{\alpha^{\prime}}\cdot P_{i}\delta_{i}^m)
=\alpha^{\prime}+r-1-m.
\end{align*}

Let $M$ be a right graded $\DS(I)$-submodule of $E_{i}$. 
Set $X_{j}:=\{\ell\mid (j,\ell)\in \Exp (M)\}$. 
From Lemma \ref{lem<exp}, there exists the smallest integer 
$j$ with $\sharp X_{j}=\infty$.
Put $s:=s_{M}:=\min\{j\mid\sharp X_{j}=\infty\}$, 
and set $M_{s}:=\{\theta\in M\mid
\exp(\theta)=(s,\ell)\ {\rm for\ some\ }\ell\}$. 
Then it is clear that $s\geq r-1$.

Let $\theta_{m}\in M$ be a homogeneous operator 
satisfying $\exp(\theta_{m})=(s,m)$ with $m\geq s$. 
Since $m-s+r-1>0$, we can write
\begin{align}
\theta_{m}=\sum_{\ell=0}^{s-r+1}
a_{\ell}y^{s-r+1-\ell}P_{i}\delta_{i}^{m-\ell}
\qquad (a_{\ell}\in K).\label{theta}
\end{align}
We may assume that $a_{0}=1$. Set 
$\Omega:=\{1,\dots,i-1,i+1,\dots,r\}$. 
For $0\leq \ell\leq s-r+1$, we write 
\begin{align*}
&\delta_{i}^{m-\ell}P_{i}\\
&=\sum^{r-1}_{\ell^{\prime}=0}[m-\ell]_{\ell^{\prime}}
\left(\frac{1}{\ell^{\prime}!(r-1-\ell^{\prime})!}
\sum_{\sigma \in S^{\Omega}}\delta(p_{\sigma(1)})
\cdots \delta(p_{\sigma(\ell^{\prime})})
p_{\sigma(\ell^{\prime}+1)}\cdots p_{\sigma(r)}\right)
\delta_{i}^{m-\ell-\ell^{\prime}}\\
&\equiv \sum^{r-1}_{\ell^{\prime}=0}[m-\ell]_{\ell^{\prime}}
d_{\ell^{\prime}}y^{r-1-\ell^{\prime}}
\delta_{i}^{m-\ell-\ell^{\prime}}
\qquad\qquad ({\rm mod}\ p_{i}D(S)),
\end{align*}
for some $d_{\ell^{\prime}}\in K$ by Lemma \ref{cal}. 
We obtain $d_{0}\neq 0,d_{r-1}\neq 0$ 
from this computation. Then
\begin{align}
\theta_{m}\cdot P_{i}\delta_{i}^{m^{\prime}}
&=\sum_{\ell=0}^{s-r+1}
a_{\ell}y^{s-r+1-\ell}P_{i}(\delta_{i}^{m-\ell}
P_{i})\delta_{i}^{m^{\prime}}\notag \\
&=\sum_{\ell=0}^{s-r+1}\sum^{r-1}_{\ell^{\prime}=0}
a_{\ell}[m-\ell]_{\ell^{\prime}}d_{\ell^{\prime}}
y^{k-\ell-\ell^{\prime}}P_{i}
\delta_{i}^{m-\ell-\ell^{\prime}}\notag \\
&=\sum_{t=0}^{s}c_{t}y^{s-t}P_{i}
\delta_{i}^{m+m^{\prime}-t},\label{pthe}
\end{align}
where
\begin{align}
c_{t}:=\sum_{{0\leq \ell\leq s-r+1 \atop 
0\leq \ell^{\prime}\leq r-1}\atop 
\ell+\ell^{\prime}=t}a_{\ell}[m-\ell]_{\ell^{\prime}}
d_{\ell^{\prime}}\label{c_t}
\end{align}
for $0\leq t\leq s$. We remark that $c_{t}$ 
does not depend on $m^{\prime}$. 
Put $m_{0}:=\max\{\ell\mid(s-1,\ell)\Exp(M)\}+s$. 
\begin{Lem}\label{rankthe}
For $1\leq j\leq r$, there exist operators 
$\theta_{m_{1}},\dots,\theta_{m_{j}}\in M_{s}$ 
such that 
\begin{align}
\rank
\begin{pmatrix}
c_{0}^{(1)}&\cdots&c_{0}^{(j)}\\
\vdots & &\vdots \\
c_{r-1}^{(1)}&\cdots&c_{r-1}^{(j)}
\end{pmatrix}
=j,\label{mat1}
\end{align}
and $m_{0}<m_{1}<\cdots<m_{j}$, 
where $c_{t}^{(j)}$ for $\theta_{m_{j}}$ 
has been defined in $(\ref{c_t})$.
\end{Lem}
\begin{proof}
We prove the assertion by induction. 
It is clear in the case $j=1$. 

Let $1<j<r$. 
Assume that there exist 
$\theta_{m_{1}},\dots,\theta_{m_{j}}\in M_{s}
\ (m_{1}<\cdots<m_{j})$ 
satisfying the condition (\ref{mat1}). 

For $m>m_{j}$, put a vector
$$
\bm{w}:=\left( y^{s}P_{i}\delta_{i}^{m+m^{\prime}},
y^{s-1}P_{i}\delta_{i}^{m+m^{\prime}-1},\dots,
P_{i}\delta_{i}^{m+m^{\prime}-s}\right),
$$
and put an $(s+1)\times(s-r+j+2)$ matrix
\begin{align*}
A:=
\scriptsize{
\begin{pmatrix}
c_{0}^{(1)}&\cdots&c_{0}^{(j)}&d_{0}&0
&0&\cdots&0 \\
c_{1}^{(1)}&\cdots&c_{1}^{(j)}&[m]_{1}d_{1}
&d_{0}&0&\cdots&0 \\
c_{2}^{(1)}&\cdots&c_{2}^{(r-1)}&[m]_{2}d_{2}
&[m-1]_{1}d_{1}&d_{0}&\ddots&0 \\
\vdots & &\vdots &\vdots& &\ddots&\ddots&\vdots \\
c_{r-2}^{(1)}&\cdots&c_{r-2}^{(j)}&[m]_{r-2}d_{r-2}
& & & &0\\
c_{r-1}^{(1)}&\cdots&c_{r-1}^{(j)}&[m]_{r-1}d_{r-1}
&[m-1]_{r-2}d_{r-2}& & &d_{0}\\
\vdots& &\vdots&0&[m-1]_{r-1}d_{r-1}& & &\vdots\\
\vdots& &\vdots&0&0&\ddots& &\vdots\\
\vdots & &\vdots &\vdots& &\ddots&\ddots&\vdots\\
c_{s}^{(1)}&\cdots&c_{s}^{(j)}&0&\cdots
&\cdots&0&[m-s+r-1]_{r-1}d_{r-1}
\end{pmatrix}
}.
\end{align*}
We consider $m$ as a variable. 
By the induction hypothesis, there exists a nonzero 
$j$-minor of the matrix in (\ref{mat1}). 
We denote by $B$ the matrix of this $j$-minor. 
We take the lowest $s-r$ rows of $A$ 
and $j$ rows from the remaining 
$r-1$ rows of $A$ so that we get the $(s-r+j+2)$-minor 
$C$ whose matrix contains the matrix $B$. 
The coefficient of the leading term of 
$C$ is the determinant of $B$. 
Thus $C$ is not zero as a polynomial in variable $m$, 
and hence the solutions of $C=0$ is finite. 
Because of this, the number of $m$ 
with $\rank(A)<s-r+j+2$ is finite. 
Hence we can take a positive integer $m>m_{j}$ 
such that $\exp(\theta_{m})\in M_{s}$, and 
$\rank(A)=s-r+j+2$. 

We write 
$\theta_{m}=\sum_{\ell=0}^{s-r+1}
a_{\ell}y^{s-r+1-\ell}P_{i}\delta_{i}^{m-\ell}$
in the same way as in (\ref{theta}). Put 
$$
\bm{v}:=\left(\lambda_{1},\dots,\lambda_{j},
\lambda_{j+1}a_{0},\dots,\lambda_{j+1}a_{s-r+1} \right).
$$
Then
\begin{align*}
\bm{w}A\,{}^t \bm{v}=
\lambda_{1}\theta_{m_{1}}\cdot P_{i}\delta_{i}^{m_{1}^{\prime}}
+\cdots+
\lambda_{j}\theta_{m_{j}}\cdot P_{i}\delta_{i}^{m_{j}^{\prime}}
+\lambda_{j+1}\theta_{m}\cdot P_{i}\delta_{i}^{m^{\prime}}
\end{align*}
with 
$m_{1}+m_{1}^{\prime}=\cdots=m_{j}+m_{j}^{\prime}=m+m^{\prime}$. 
If $\bm{w}A\,{}^t \bm{v}=0$, then $\bm{v}=0$ 
since $\rank(A)=s-r+j+2$. Therefore 
$\{ \theta_{m_{1}}\cdot P_{i}\delta_{i}^{m_{1}^{\prime}},
\dots,\theta_{m_{j}}\cdot P_{i}\delta_{i}^{m_{j}^{\prime}},
\theta_{m}\cdot P_{i}\delta_{i}^{m^{\prime}}\}$ 
is linearly independent over $K$. 

Put $\theta_{m_{j+1}}:=\theta_{m}$, and 
suppose that 
\begin{align*}
\rank
\begin{pmatrix}
c_{0}^{(1)}&\cdots&c_{0}^{(j+1)}\\
\vdots & &\vdots \\
c_{r-1}^{(1)}&\cdots&c_{r-1}^{(j+1)}
\end{pmatrix}
<j+1.
\end{align*}
Then there exists $(\lambda_{1},\dots,\lambda_{j+1})
\in K^{j+1}\setminus\{\bm{0}\}$ such that 
\begin{align}
\lambda_{1}\,{}^t(c_{0}^{(1)},\dots,c_{r-1}^{(1)})
+\cdots+
\lambda_{j+1}\,{}^t(c_{0}^{(j+1)},\dots,c_{r-1}^{(j+1)})
=\bm{0}.\label{l-i1}
\end{align}
Since $\{ \theta_{m_{1}}\cdot P_{i}\delta_{i}^{m_{1}^{\prime}},
\dots,\theta_{m_{j+1}}\cdot P_{i}\delta_{i}^{m_{j+1}^{\prime}}\}$ 
is linearly independent, we have 
\begin{align}
\sum_{k=0}^{j+1}\lambda_{k}
\theta_{m_{k}}\cdot P_{i}\delta_{i}^{m_{k}^{\prime}}\neq 0.
\label{l-i2}
\end{align}
Hence we can write 
$\exp\left(\sum_{k=0}^{j+1}\lambda_{k}
\theta_{m_{k}}\cdot P_{i}\delta_{i}^{m_{k}^{\prime}}
\right)=(\alpha,\beta)$ 
for some $\alpha<s$ and 
$\beta>m_{j+1}+m_{j+1}^{\prime}-s
>m_{0}-s=\max\{\ell\mid(s-1,\ell)\Exp(M)\}$ 
by (\ref{l-i1}) and (\ref{l-i2}). 
This is a contradiction.
\end{proof}
Let $M$ be a right graded $\DS(I)$-submodule 
of $E_{i}$. For a nonnegative integer $\ell$ with 
$(k,\ell)\in\Exp(M)$ for some $k$, 
we define an integer $t_{\ell}$ by
$$
t_{\ell}:=\min\{k\mid(k,\ell)\in\Exp(M)\}.
$$

By Lemma \ref{rankthe}, there exist operators 
$\theta_{m_{1}},\dots,\theta_{m_{r}}\in M_{s}$ 
satisfying the condition (\ref{mat1}). 
We denote by $N$ the right submodule of $M$ 
generated by the operators 
$\theta_{m_{1}},\dots,\theta_{m_{r}}$. 
\begin{Lem}\label{lemmain}
There exists a positive integer $n_{0}$ such that, 
for any $m\geq n_{0}$, 
\begin{align*}
(s,m)\in\Exp(N)\ {\rm and}\ t_{m}=s.
\end{align*}
\end{Lem}
\begin{proof}
By Lemma \ref{rankthe}, 
there exist 
$\theta_{m_{1}},\dots,\theta_{m_{r}}\in M_{s}$ 
such that 
\begin{align*}
\rank
\begin{pmatrix}
c_{0}^{(1)}&\cdots&c_{0}^{(r)}\\
\vdots & &\vdots \\
c_{r-1}^{(1)}&\cdots&c_{r-1}^{(r)}
\end{pmatrix}
=r,\ \ {\rm and}\ \ \rank
\begin{pmatrix}
c_{0}^{(1)}&\cdots&c_{0}^{(r)}\\
\vdots & &\vdots \\
c_{r-2}^{(1)}&\cdots&c_{r-2}^{(r)}
\end{pmatrix}
<r.
\end{align*}
Then there exists a nonzero vector 
$(\lambda_{1},\dots,\lambda_{r})
\in K^{r}\setminus\{\bm{0}\}$ 
such that 
$$
\lambda_{1}\,{}^t(c_{0}^{(1)},\dots,c_{r-2}^{(1)})
+\cdots+
\lambda_{r}\,{}^t(c_{0}^{(r)},\dots,c_{r-2}^{(r)})
=\bm{0},
$$
and
$$
\lambda_{1}c_{r-1}^{(1)}+\cdots+
\lambda_{r}c_{r-1}^{(r)}\neq 0.
$$
Put $\theta:=\sum_{k=0}^{r}\lambda_{k}
\theta_{m_{k}}\cdot P_{i}\delta_{i}^{m_{k}^{\prime}}
\in N$ with 
$m_{1}+m_{1}^{\prime}=\cdots=m_{r}+m_{r}^{\prime}$. 
It follows that 
$\exp(\theta)=(s,m_{r}+m_{r}^{\prime}-r+1)$. 
Put $n_{0}=m_{r}-r+1$, and 
put $m_{r}^{\prime}=m-n_{0}$ for any $m\geq n_{0}$. 
Thus 
$$
(s,m)=(s,m_{r}+m_{r}^{\prime}-r+1)=
\exp(\theta)\in\Exp(N).
$$

It remains to prove that $t_{m}=s$. 
We have $t_{m}\geq s$ since 
$m\geq n_{0}\geq m_{0}-s=\max\{\ell\mid(s-1,\ell)\Exp(M)\}$. 
Conversely we have 
$t_{m}\leq s$ since $(s,m)\in\Exp(M)$, as required.
\end{proof}
Let $R$ be a graded ring. A right graded $R$-module 
$M$ is said to be right gr-Noetherian, if $M$ satisfies 
the ascending chain condition for graded submodules 
of $M$. It is straightforward to verify that 
$M$ is right gr-Noetherian if and only if 
each graded submodule of $M$ is finitely generated. 
\begin{Prop}
The right $\DS(I)$-module $E_{i}$ is right Noetherian.
\end{Prop}
\begin{proof}
Recall that $\DS(I)$ is a graded ring 
by the total degree, and that 
$E_{i}$ is a graded $\DS(I)$-module. 
By \cite[Theorem I\hspace{-.1em}I.3.5]{Nast}, 
it is enough to prove that $E_{i}$ is right 
gr-Noetherian. Let $M$ be a right graded $\DS(I)$-submodule 
of $E_{i}$. We will prove that 
$M$ is finitely generated. 

Set 
$$
G:=\left\{ (t_{\ell},\ell)\mid\ell<n_{0}
\ {\rm and}\ (k,\ell)\in\Exp(M)\ 
{\rm\ for\ some}\ k\right\}.
$$
Then $G$ is a finite set. Fix an operator 
$\theta_{(t_{\ell},\ell)}\in M$ for 
$(t_{\ell},\ell)\in G$, and set
$$
\overline{G}:=\left\{\theta_{(t_{\ell},\ell)}
\in M\mid (t_{\ell},\ell)\in G\right\}.
$$
Then $\overline{G}$ is also a finite set. 

Let $n_{0}$ be the integer satsfying Lemma \ref{lemmain}. 
We denote by $M^{\prime}$ the right $\DS(I)$-module 
generated by $\overline{G}$ and $N$. 
Then $M^{\prime}$ is finitely generated and 
$M^{\prime}\subseteq M$. 

Let $(k,m)\in\Exp(M)$, then $k\geq t_{m}$. 
If $m<n_{0}$, then $(t_{m},m)\in G\subseteq\Exp(M^{\prime})$ 
by the definitions of $t_{m}$ and $G$. We have 
$(k,m)=(t_{m}+k-t_{m},m)\in\Exp(M^{\prime})$
by Lemma \ref{lem<exp}. 

If $m\geq n_{0}$, then $(s,m)\in\Exp(M^{\prime})$ by 
Lemma \ref{lemmain}. 
It follows from Lemma \ref{lem<exp} that 
$(k,m)=(s+k-s,m)\in\Exp(M^{\prime}).$
Hence $\Exp(M^{\prime})=\Exp(M)$. 
The assertion follows from Lemma \ref{lemexp}.
\end{proof}
\begin{Cor}\label{noeth>r-1}
The right $\DS(I)$-module 
$$
\left(L_{i-1}/L_{i}\right)^{\geq r-1}
\cdot \DS(I)=\left(L_{i-1}^{\geq r-1}
\cdot\DS(I)+L_{i}/L_{i}\right)
$$ 
is right Noetherian.
\end{Cor}
Next we study the $S$-module 
$\big{(}L_{i-1}/L_{i}\big{)}^{< r-1}$. 
\begin{Lem}\label{Li<r-1}
The $K$-vector space
\begin{align*}
L_{i}^{<r-1}:=\bigoplus_{m<r-1}L_{i}^{(m)}
\end{align*}
is a right $S$-module.
\end{Lem}
\begin{proof}
Suppose that $0\leq m<r-1$. Let 
$\theta \in L_{i}^{(m)}\subseteq\DS(I)$. 
For $f\in S$, 
$$
\theta f(QS)=\theta(QfS)\subseteq I.
$$
Thus $\theta f\in\DS(I)$. 
It follows from Proposition \ref{p1} that 
$\theta f\in \bigoplus_{\ell=0}^{m}\DS^{(\ell)}(I)$. 
The operator $\theta f$ is divisible by 
the polynomial $p_{1}\cdots p_{i}$ since 
$\theta\in p_{1}\cdots p_{i}\DS^{(m)}(S)$. 
Thus each homogeneous component of $\theta f$ 
is divisible by $p_{1}\cdots p_{i}$. 
It follows that 
$$
\theta f\in\bigoplus_{\ell=0}^{m}
\big{(}\DS^{(\ell)}(I)\cap
(p_{1}\cdots p_{i})\DS^{(\ell)}(S)\big{)}
=\bigoplus_{\ell=0}^{m}L_{i}^{(\ell)}.
$$
Hence $L_{i}^{<r-1}\cdot S\subseteq L_{i}^{<r-1}$.
\end{proof}
The following holds in general.
\begin{Prop}\label{Weylp}
As a vector space, 
$$
\bigoplus_{|\bm{\alpha}|<r-1}S\partial^{\bm{\alpha}}
=\bigoplus_{|\bm{\alpha}|<r-1}\partial^{\bm{\alpha}}S.
$$
\end{Prop}
Define a right $S$-module 
$\DS(S)^{<r-1}:=\bigoplus_{|\bm{\alpha}|<r-1}
\partial^{\bm{\alpha}}S$. Then $\DS(S)^{<r-1}$ is the 
module of differential operators 
of order less than $r-1$ by Proposition \ref{Weylp}. 
By Lemma \ref{Li<r-1}, we have 
the inclusion of right $S$-modules: 
\begin{align*}
L_{i}^{<r-1}\subseteq\DS(S)^{<r-1}.
\end{align*}
\begin{Lem}\label{S-noeth}
The right $S$-module $\big{(}L_{i-1}/L_{i}\big{)}^{<r-1}$ 
is Noetherian.
\end{Lem}
\begin{proof}
Since $\DS(S)^{<r-1}$ is 
a finitely generated right $S$-module, 
$\DS(S)^{<r-1}$ is Noetherian as a right $S$-module. 
Hence the subquotient 
$\big{(}L_{i-1}/L_{i}\big{)}^{<r-1}=
\big{(}L_{i-1}^{<r-1}+L_{i}/L_{i}\big{)}$ 
of $\DS(S)^{<r-1}$ 
is Noetherian as a right $S$-module.
\end{proof}
\begin{Lem}\label{DS-noeth}
The right $\DS(I)$-module $L_{i-1}/L_{i}$ 
is Noetherian.
\end{Lem}
\begin{proof}
Let $M$ be a right $\DS(I)$-submodule of $L_{i-1}/L_{i}$. 
Put $M^{<r-1}:=M\cap \big{(}L_{i-1}^{<r-1}+
L_{i}/L_{i}\big{)}$, and put 
$M^{\geq r-1}:=M\cap \big{(}L_{i-1}^{\geq r-1}
+L_{i}/L_{i}\big{)}$. 
Then $M^{<r-1}$ is a right $S$-submodule of 
$\big{(}L_{i-1}^{<r-1}+L_{i}/L_{i}\big{)}$. 
By Lemma \ref{S-noeth}, $M^{<r-1}$ is finitely generated 
as a right $S$-module. Then 
$M^{<r-1}\cdot \DS(I)=M\cap
\big{(}L_{i-1}^{<r-1}\cdot\DS(I)+L_{i}/L_{i}\big{)}$ 
is also finitely generated as a right $\DS(I)$-module 
with generators $u_{1},\dots ,u_{\ell}$. 

As right $\DS(I)$-modules,
\begin{align*}
\big{(}M/M^{<r-1}\cdot\DS(I)\big{)}&\simeq 
\frac{M+\big{(}L_{i-1}^{<r-1}\cdot\DS(I)
+L_{i}/L_{i}\big{)}}{\big{(}L_{i-1}^{<r-1}
\cdot\DS(I)+L_{i}/L_{i}\big{)}}\\
&\subseteq\frac{\big{(}L_{i-1}/L_{i}\big{)}}
{\big{(}L_{i-1}^{<r-1}\cdot\DS(I)
+L_{i}/L_{i}\big{)}}\\
&\simeq \frac{\big{(}L_{i-1}^{\geq r-1}
\cdot\DS(I)+L_{i}/L_{i}\big{)}}
{\big{(}L_{i-1}^{<r-1}\cdot\DS(I)+L_{i}
/L_{i}\big{)}\cap\big{(}L_{i-1}^{\geq r-1}
\cdot\DS(I)+L_{i}/L_{i}\big{)}}.
\end{align*}
By Corollary \ref{noeth>r-1}, 
$\big{(}L_{i-1}^{\geq r-1}\cdot\DS(I)
+L_{i}/L_{i}\big{)}$ is Noetherian as a 
right $\DS(I)$-module. Thus the module 
$\big{(}M/M^{<r-1}\cdot\DS(I)\big{)}$ 
is finitely generated as a $\DS(I)$-module. 
We take a system of generators 
$\{\overline{v_{1}},\dots,
\overline{v_{\ell^{\prime}}}\}$.

We will prove that $M$ is 
a finitely generated $\DS(I)$-module.
Fix $\theta \in M$. Since $\overline{\theta}\in
\big{(}M/M^{<r-1}\cdot \DS(I)\big{)}$, 
we can write 
$$
\overline{\theta}=\sum_{j}
\overline{v_{j}}\cdot\eta_{j}
$$
for some $\eta_{j}\in\DS(I)$. 
We also write 
$$
\theta-\sum_{j}v_{j}\cdot\eta_{j}=
\sum_{k}u_{k}\cdot\eta^{\prime}_{k}
$$
for some $\eta^{\prime}_{k}\in\DS(I)$, 
since $\theta-\sum_{j}v_{j}\eta_{j}
\in M^{<r-1}\cdot \DS(I)$. Thus 
$$
\theta\in\sum_{j}v_{j}\cdot\DS(I)
+\sum_{k}u_{k}\cdot\DS(I).
$$
Therefore $M=\sum_{j}v_{j}\cdot\DS(I)
+\sum_{k}u_{k}\cdot\DS(I)$ is finitely generated.
\end{proof}
\begin{Them}
The ring $\DS(S/I)\simeq\DS(I)/I\DS(S)$ 
is Noetherian (i.e., $\DS(S/I)$ is 
right Noetherian and left Noetherian).
\end{Them}
\begin{proof}
By Lemma \ref{DS-noeth} and by 
considering the sequence (\ref{seqL}), 
we see that the ring $\DS(I)/I\DS(S)$ is 
right Noetherian. Therefore, by Corollary \ref{r-l-n}, 
we conclude that the ring 
$\DS(S/I)\simeq\DS(I)/I\DS(S)$ is Noetherian.
\end{proof}
In the rest of this section, we give an 
example of a family of Noetherian rings whose 
graded rings associated to the order filtration 
are not Noetherian. 

By Proposition \ref{p1}, we can decompose 
$\DS(I)/I\DS(S)$ into the direct sum 
$$
\DS(I)/I\DS(S)=\bigoplus_{m\geq 0}\big{(}
\DS^{(m)}(I)/I\DS^{(m)}(S)\big{)}
$$
as a left $S$-module. The order filtration of 
$\DS(I)/I\DS(S)$ is the filtration 
$\FC=\{F_{m}\}_{m\geq 0}$ defined by 
$$
F_{m}=\bigoplus_{\ell\leq m}\big{(}
\DS^{(\ell)}(I)/I\DS^{(\ell)}(S)\big{)}.
$$

We denote by $S_{j}$ the $K$-vector subspace of $S$ 
spanned by the monomials of degree $j$. 
An element 
$\theta=\sum_{\bm{\alpha}}f_{\bm{\alpha}}
\partial^{\bm{\alpha}}\in\DS(S)$ 
is of polynomial degree $k$, 
if $k$ is the smallest integer such that 
$f_{\bm{\alpha}}\in\bigoplus_{j=0}^{k}S_{j}$ 
for all $\bm{\alpha}$ with nonzero $f_{\bm{\alpha}}$.
\begin{Examp}\label{example}
Let $S=k[x,y]$ be the polynomial ring, 
and let $I$ be the ideal generated by 
the polynomial $Q=p_{1}\cdots p_{r}\ (r\geq 2)$
 defining a central arrangement. 

The graded ring $\Gr\DS(S/I)$ associated to 
the order filtration is a commutative ring. 
We consider the ideal 
$M:=\langle \overline{P_{1}\delta_{1}^m}
\mid m\geq 1\rangle$ 
of $\Gr\DS(S/I)$. 

Assume that $M$ is finitely generated with 
generators $\eta_{1},\dots,\eta_{\ell}$. 
Then there exists a positive integer $m$ such that 
$$
M=\langle\eta_{1},\dots,\eta_{\ell}\rangle
\subseteq\langle\overline{P_{1}\delta_{1}}
,\dots,\overline{P_{1}\delta_{1}^{m-1}}\rangle.
$$
Since $\overline{P_{1}\delta_{1}^{m}}\in M$, 
we can write
\begin{align}
\overline{P_{1}\delta_{1}^m}=
\overline{P_{1}\delta_{1}}
\cdot\overline{\theta_{1}}+
\cdots+\overline{P_{1}\delta_{1}^{m-1}}
\cdot\overline{\theta_{m-1}}\label{Pd^m}
\end{align}
for some 
$\theta_{1},\dots,\theta_{m-1}\in\DS(I)$.

If $\theta\in\DS(I)$ with $\ord(\theta)\leq 1$, 
then the polynomial degree of $\theta$ is 
greater than or equal to $1$ by Proposition 
\ref{pb-2}. 
Since the order of the LHS of (\ref{Pd^m}) 
equals $m$, there exists at least one 
$\theta_{j}$ such that the order of 
$\theta_{j}$ is greater than or equal to $1$. 
Thus the polynomial degree of 
the RHS of (\ref{Pd^m}) is greater than $r-1$. 
However, the polynomial degree of the LHS 
of (\ref{Pd^m}) is exactly $r-1$. 
This is a contradiction. 

Therefore $M$ is not finitely generated, 
and thus we have proved that $\Gr\DS(S/I)$ 
is not Noetherian. 
\end{Examp}

\vspace{10mm}
\end{document}